\documentclass[12pt]{article}
\usepackage{amssymb,latexsym,amsmath}
\usepackage{graphics}                 
\usepackage{color}                    
\usepackage{hyperref}                 
\usepackage[all]{xy}

\begin{document}
\pagestyle{plain} 

\newtheorem{theorem}{Theorem}[section]
\newtheorem{lemma}[theorem]{Lemma}
\newtheorem{proposition}[theorem]{Proposition}
\newtheorem{corollary}[theorem]{Corollary}

\newtheorem{definition}[theorem]{Definition}
\newtheorem{example}[theorem]{Example}
\newtheorem{xca}[theorem]{Exercise}
\newtheorem{conjecture}[theorem]{Conjecture}

\newtheorem{remark}[theorem]{Remark}

\newtheorem{correction}{Correction}

\def\bp{{\bf Proof.}\hspace{2mm}}
\def\qe{\hfill$\Box$}

\def\A{{\mathbb{A}}}
\def\C{{\mathbb{C}}}
\def\L{{\mathbb{L}}}
\def\O{{\mathcal{O}}}
\def\P{{\mathbb{P}}}
\def\Q{{\mathbb{Q}}}
\def\Z{{\mathbb{Z}}}
\def\cZ{{\mathcal{Z}}}
\def\Ch{{\rm Ch}}
\def\p{{\mathbf{p}}}

\def\sp{{\rm Sp}}
\def\rank{{\rm rank}}
\def\HC{{\rm Hodge}}
\def\GHC{{\rm GHC}}
\def\Griff{{\rm Griff}}

\title{The Generalized Bloch Conjecture for the quotient of certain Calabi-Yau varieties}
\author{Wenchuan Hu}
\maketitle

\begin{abstract}
In this paper,  the generalized  Bloch Conjecture on zero cycles for the
quotient of certain complete intersections with trivial canonical bundle is proved to hold.

As an application of Bloch-Srinivas method on the decomposition of the diagonal,
we compute the rational coefficient Lawson homology for 1-cycles
and codimension two cycles for these  quotient varieties. The (Generalized) Hodge Conjecture
is proved to hold for codimension two cycles (and hence also for 2-cycles) on these quotient
varieties.
\end{abstract}

\tableofcontents

\section{Introduction}

In this paper, all varieties are defined over $\C$.
For a projective variety $X$, denote by  $\cZ_p(X)$ the spaces of algebraic $p$-cycles
and $\Ch_p(X)$  the Chow group of
$p$-cycles on $X$, i.e, $\Ch_p(X)=\cZ_p(X)/{\hbox{\{rational equivalence\}}}$. Let $cl_p:\Ch_p(X)\to H_{2p}(X,\Z)$ be the cycle class map.
Tensoring with ${{\Q}}$, we have
$cl_p\otimes{\Q}:{\Ch}_p(X)\otimes{\Q}\to H_{2p}(X,\Q).$
Let $\Ch_p(X)_{hom}\subset \Ch_p(X)$ be the subgroup of $p$-cycles homologous to zero. Set $\Ch^q(X):=\Ch_{n-q}(X)$.

In 1968, D. Mumford showed that $\Ch_0(X)_{hom}$ is not finite dimensional
for a smooth projective surface $X$ with non-vanishing geometric genus $p_g(X)$ (cf. \cite{Mumford}).
This result was generalized by Ro\u\i tman to arbitrary dimension (cf. \cite{Roitman}). In this situation,
a nontrivial conjecture of Bloch   asserts  that  if  a smooth projective surface $X$ with $p_g(X)=0$, then
$\Ch_0(X)_{hom}$ is finite dimensional (\cite{Bloch}). Equivalently, if $p_g(X)=0$, then there is a curve $C\subset X$ such
that the natural map $\Ch_0(C)\to \Ch_0(X)$ is surjective.

 This conjecture can be generalized as follows (cf. \cite{Paranjape-Srinivas}):

\begin{conjecture}[Generalized Bloch Conjecture]
Let $X$ be a  smooth projective variety satisfying $H^{i,0}(X) = 0$ for all $i>r$.
Then there is a subvariety $i : Z \hookrightarrow  X$, where $\dim Z = r$, such that $i_* : \Ch_0(Z) \to \Ch_0(X)$
is surjective.
\end{conjecture}

Some examples are known in support of these conjectures. For example, Bloch's conjecture is true for surfaces which are not of general type \cite{Bloch-K-L}.
This conjecture also holds for some surfaces  of general type which are quotients of some special surfaces by a free finite group action
(cf. \cite{Inose-Mizukami}, \cite{Voisin}). In higher dimensional case, it was proved by Ro\u\i tman \cite{Roitman} that $\Ch_0(X)\cong \Z$
for smooth projective complete intersection with $H^{i,0}(X) = 0$ for all $i>0$. In \cite{Bloch-S}, it was proved by Bloch and Srinivas that the
Generalized Bloch Conjecture holds for  Kummer varieties of odd dimensions.

In this paper, the Generalized Bloch Conjecture is proved  to hold for  quotients of certain even dimensional complete intersection with
trivial canonical bundle by a free involution  and for the resolution of singularities
to the quotients of certain odd dimensional complete intersection with trivial canonical bundle by an involution with isolated fixed points.
Based on these results, we compute the rational Lawson homology and verify the Generalized Hodge
Conjecture  on 1-cycles and codimension-2  cycles   for these varieties.

\medskip
I would like to express my gratitude to Michael Artin for helpful discussion
and suggestion during the preparation of this paper.

\section{Main results}
Now let $X^n\subset\P^{2n+1}$ be the complete intersection of quadrics
$$Q_i(z_0,z_1\cdots,z_n)+Q_i'(z_{n+1},z_{n+2}\cdots,z_{2n+1})=0, i=0,1,\cdots,n, $$
where $Q_i$ and $Q_i'$ are quadratic forms in $n+1$ variables. For our propose, we assume
\begin{equation}\label{eq1}
Q_i(z_0,z_1\cdots,z_n)=\sum_{j=0}^na_{ij}z_j^2
\end{equation}
for $i=0,1,\cdots,n$. We also assume that $X^n$ is smooth, which holds for
the generic choice of $Q_i'$ and the choice of $a_{ij}, i,j=0,1,\cdots, n$
such that $\det(a_{ij})\neq 0$. From the direct calculation we know $X^n$ is a
Calabi-Yau $n$-fold, i.e., $K_{X^n}$ is trivial and hence
$h^{n,0}(X^n):=\dim H^{n,0}(X^n)=1$.

For $n=2m$ a positive even integer, we define an involution of $\P^{2n+1}$ by
 \begin{equation}\label{eq2}
 \sigma:(z_0:z_1:\cdots:z_{2n+1})\mapsto (-z_0: -z_1:\cdots: -z_n, z_{n+1},\cdots, z_{2n+1})
\end{equation}
which takes $X^n$ to itself. The quotient $Y^n=X^n/\langle \sigma\rangle$ is
a smooth projective variety with $H^{i,0}(Y^n)=0$ for all $i\geq 1$ for
$n=2m$  even (cf. Lemma \ref{lemma2.2}).  Denote by $\pi:X^n\to Y^n$  the projection.

Our first main result is following theorem.
\begin{theorem}\label{Th1.1}
Let $n$ be a positive even integer. The Generalized Bloch Conjecture holds for $Y^n$, i.e.,
for the projective variety $Y^n=X^n/\langle  \sigma\rangle$ above, we have $\Ch_0(Y^n)=\Z$.
\end{theorem}

If $n=2m-1$ is a positive odd integer, then
we define another involution of $\P^{2n+1}$ by
 \begin{equation}\label{eq3}
 \rho:(z_0:z_1:\cdots:z_{2n+1})\mapsto (z_0: -z_1:\cdots: -z_n, z_{n+1},\cdots, z_{2n+1})
\end{equation}
  which takes
 $X^n$ to itself. The involution $\rho$ has $2^{2m}$  isolated fixed points.
 Hence the quotient $Y^{2m-1}=X^{2m-1}/\langle \rho\rangle$ is a  projective variety with $2^{2m}$ isolated singular points, denote by $q_i, i=1,2,\cdots,2^{2m}$.  Denote also by $\pi:X^{2m-1}\to Y^{2m-1}$  the projection. Each singular point is a cyclic quotient singular point.
 Let $\widetilde{Y^{2m-1}}\to Y^{2m-1}$ be a resolution of singularity, then
the exceptional divisor $E_i$ at each singular point $q_i$ has
 only normal crossings in
$\widetilde{Y^{2m-1}}$ and every irreducible component of $E_i$ is nonsingular and rational \cite{Fujiki}.

Our second main result is the following theorem.
\begin{theorem}\label{Th1.7}
Let $n=2m-1$ be a positive odd integer.
The Generalized Bloch Conjecture holds for $\widetilde{Y^{n}}$, i.e., for $\widetilde{Y^{n}}$ above,
we have  $\Ch_0(\widetilde{Y^{2m-1}})\cong\Z$.
Moreover,  $\Ch_0(Y^{2m-1})\cong\Z$.
\end{theorem}

The application of the main results on algebraic cycles and Lawson homology is given in section \ref{sec4}.

\section{The proof of  main theorems}\label{sec3}
\begin{lemma}
For the generic choice of $Q_i'$ and the choice of $a_{ij}, i,j=0,1,\cdots, n$ such that
the determinant $\det(a_{ij})$ of the matrix $(a_{ij})$ is nonzero, then
$X^{n}$ is a smooth projective variety of dimension $n$.
\end{lemma}
\bp It follows from the definition of smoothness of projective variety.
Note that $\det(a_{ij})\neq 0$ implies that there is no common solution
 for the system of equations by  $Q_i=0$ and those of partial derivatives.
 Similarly for the generic choice of $Q_i'$, there is no common solution for the system of equations $Q_i', i=1,2,\cdots, n+1$.
\qe

\begin{lemma}\label{lemma2.2}
For $n=2m$ a positive even integer, the quotient $Y^n=X^n/\langle \sigma\rangle$ by the involution $\sigma:X^n\to X^n$
is a smooth projective variety. Moreover, $H^{i,0}(Y^n)=0$ for all $i>0$.
\end{lemma}
\bp
The involution $\sigma:X^n\to X^n$ is induced by the involution of
$\P^{2n+1}$ defined by $\sigma:(z_0:z_1:\cdots:z_{2n+1})\mapsto (-z_0: -z_1:\cdots: -z_n, z_{n+1},\cdots, z_{2n+1})$.
 By the assumption, the fixed point set  of $\sigma:\P^{2n+1}\to \P^{2n+1}$
 and $X^n$ have no intersection since the system of equations $Q_i=0$, $i=0,1,\cdots, n$
has no common solution in $\P^n$ by the assumption that $\det(a_{ij})\neq 0$. Similarly for a
generic choice of $Q_i'$, $i=0,1,\cdots, n$. Therefore,
$\sigma:\P^{2n+1}\to \P^{2n+1}$
induces a fixed point free involution on $X$ and hence the quotient
$Y^n=X^n/\langle \sigma\rangle$ is a smooth projective variety.
Note that $H^{i,0}(Y^n)=0$ for $i> n$ for the reason of dimension.
For $i<n$, $\dim H^{i,0}(Y^n)\leq \dim H^{i,0}(X^n)$ and the latter is zero by
Lefschetz hyperplane Theorem. For $i=n$, we have
$2=\chi(\O_{X^n})=2 \chi (\O_{Y^n})=2(1-\dim H^{1,0}(X^n)+\cdots+(-1)^n\dim H^{n,0}(Y^n))=2(1+\dim H^{n,0}(Y^n))$, where the second equality holds since
$X^n\to Y^n$ is a \'{e}tale morphism (cf. Example 18.3.9 in \cite{Fulton})
and the last equality holds since $n$ is an even integer and so $\dim H^{n,0}(Y^n)=0$.
\qe

\begin{remark}
Since $H^{i,0}(Y^n)=0$ for all $i>0$, the Generalized Bloch Conjecture implies
$\Ch_0(Y^n)\otimes\Q=\Q$. The statement in Theorem \ref{Th1.1} is slightly
stronger than this.
\end{remark}

A well-known result is needed in our computation.
\begin{lemma}\label{lemma2.3}
Suppose  a finite group $G$ acts on a variety $X$ with nonsingular quotient
variety $Y=X/G$. Let $\pi:X\to Y$ be the quotient map. Then there exist two homomorphisms
$\pi_*:\Ch_0(X)_{hom}\to \Ch_0(Y)_{hom}$ and $\pi^*:\Ch_0(Y)_{hom}\to \Ch_0(X)_{hom}$
such that
\begin{equation}\label{eq4}
\left\{\begin{array}{lll}
\pi^*\pi_*&=& \sum_{g\in G} g_*\\
\pi_*\pi^*&=& N
\end{array}\right.
\end{equation}
where $N=|G|$ means the multiplication by $N$ in $\Ch_0(X)_{hom}$. In particular,
$ \Ch_0(Y)_{hom}=0$ if and only if  $\sum_{g\in G} g_*=0$ in $End(\Ch_0(Y)_{hom})$.
\end{lemma}
\bp See, e.g. \cite{Inose-Mizukami}, Lemma 1 and Lemma 2,
where $X$ is a surface. The point is that both $\pi_*$ and $\pi^*$ are well-defined.  The proof  works in higher dimensional case and
the case that  $Y$ is singular (cf. \cite{Fulton}).

\qe

Now we can apply Lemma \ref{lemma2.3} to the quotient map $\pi:X^n\to Y^n$.
Let $\pi_*:\Ch_0(X^n)\to\Ch_0(Y^n)$ be the push forward map and
let $\pi^*:\Ch_0(Y^n)\to\Ch_0(X^n)$ be the pull back. Then we have
$$
\pi_*\pi^*=2:\Ch_0(Y)_{hom}\to \Ch_0(Y)_{hom}
$$
and
$$
\pi^*\pi_*=\sigma_*+1:\Ch_0(X)_{hom}\to \Ch_0(X)_{hom}
$$

Since $\Ch_0(Y^n)_{hom}$ is divisible (cf. \cite{Bloch2}, \cite{Roitman}), it suffices to show
$$2\Ch_0(Y)_{hom} = \pi_*\pi^*\Ch_0(Y)_{hom} = 0.$$ Since $\pi_*$ is surjective,
it suffices to show $\pi^*\pi_*:\Ch_0(X^n)_{hom}\to \Ch_0(X^n)_{hom}$
is the zero map. That is, we need to show that $\sigma_*=-1:\Ch_0(X^n)_{hom}\to \Ch_0(X^n)_{hom}$.

Therefore,  Theorem \ref{Th1.1} follows from the following proposition.

\begin{proposition}\label{prop2.1}
Let $\sigma:X^n\to X^n$ be induced by the involution in Equation (\ref{eq2}). Then
$\sigma_*=-1:\Ch_0(X^n)_{hom}\to \Ch_0(X^n)_{hom}$.
\end{proposition}

To prove this  proposition, we need some auxiliary results.
Let $\tau_i:X^n\to X^n$ be the automorphism of $X^n$ induced by
$\sigma_i:\P^{2n+1}\to \P^{2n+1}$, where
$$\sigma_i:(z_0:\cdots: z_{i-1}:z_i:z_{i+1}\cdots:z_{2n+1})
\mapsto (z_0:\cdots: z_{i-1}:-z_i:z_{i+1}:\cdots:  z_{2n+1})
$$
for $i=0,1,\cdots, n$. Note that $\tau_i$ maps $X^n$ to itself by the assumption of $Q_i$.
Then $\tau_i:X^n\to X^n$ induces a homomorphism on Chow group of 0-cycles $\tau_{i*}:\Ch_0(X^n)\to \Ch_0(X^n)$ and so
$\tau_{i*}:\Ch_0(X^n)_{hom}\to \Ch_0(X^n)_{hom}$. Then we have the following result.

\begin{lemma}\label{lemma2.6}
The homomorphism $\tau_{i*}=-1:\Ch_0(X^n)_{hom}\to \Ch_0(X^n)_{hom}$ for all $i=0,1,\cdots, n$.
\end{lemma}
\bp
By the symmetry of $t_0,t_1,\cdots,t_n$, we only need to show the case for $i=0$.
From the definition of $X^n$ and the assumption that the matrix
$(a_{ij})$ is non-degenerated, we can make a linear transformation such that
there is only one quadratic, say $Q_1+Q_1'$, depending on the variable $t_0$.
Then one can see that the function field of  $X^n/\langle\tau_1\rangle$
is
$$\C(t_1,t_2,\cdots,t_{2n},t_{2n+1})/\langle Q_i+Q_i'=0, i=2,3,\cdots,n+1\rangle.$$

Note that the variety $Y_0\subset \P^{2n}$ defined by
equations $Q_i+Q_i'=0,i=2,3,\cdots,n+1$ is a smooth complete intersection.
Since the sum of the degrees of the defining equations of $Y_0$ is
 $\sum_{i=1}^n 2=2n$, $Y_0$ is a smooth Fano variety. This
implies $Y_0$ is rationally connected (cf. \cite{Campana}, \cite{Kollar-Miyaoka-Mori}). Therefore, for
any two generic points
$p_1,p_2$  on $Y_0$, there is a rational curve $C$ passing through $p,q$.

From the definition, the rational function field of $Y_0$ is also
$$\C(t_1,t_2,\cdots,t_{2n},t_{2n+1})/\langle Q_i+Q_i'=0, i=2,3,\cdots,n+1\rangle.$$
So $Y_0$ is birational equivalent to $X^n/\langle\tau_1\rangle$.
Hence, for any generic two points on $X^n/\langle\tau_1\rangle$, there also exists
a curve passing through the two points, i.e.,  $X^n/\langle\tau_1\rangle$ is a rationally connected variety. So $\Ch_0 (X^n/\langle\tau_1\rangle) \cong
\Z$ and $\Ch_0(X^n/\langle\tau_1\rangle)_{hom} = 0$.

This together with Lemma \ref{lemma2.3}  implies  that $\tau_{1*}+1=0 \in End(\Ch_0(X^n))_{hom}$.
This completes the proof of the Lemma.
\qe

Note that from the definition we have $\sigma=\tau_0\circ \tau_1\circ\cdots\circ\tau_n$ and so
$\sigma_*=\tau_{0*}\circ \tau_{1*}\circ\cdots\circ\tau_{n*} = (-1)^{n+1} = -1$ since $n=2m$ is an even integer.
This completes the proof of Proposition \ref{prop2.1} and hence Theorem \ref{Th1.1}.\qe

\medskip
In the following, we focus on the proof of Theorem \ref{Th1.7}. Note that
$n=2m-1$ in the below of this section.
\begin{lemma}
For the generic choice of $Q_i'$ and $(a_{ij})$, the variety $Y^{2m-1}=X^{2m-1}/\langle\rho\rangle$ has exact $2^{2m}$ isolated singular points.
\end{lemma}
\bp
Note that the set of singular points on $Y^{2m-1}$ is exact the set of fixed
points of the involution $\rho: X^{2m-1}\to X^{2m-1}$. This  fixed points set is defined by  equations
\begin{equation}\label{eq5}
t_1=t_2=\cdots=t_n=0 \quad and \quad Q_i+Q_i'=0, ~ i=0,1,\cdots, n,
\end{equation}
i.e., the intersection of $\P^{n+1}$ and $X^n$ in $\P^{2n+1}$.
Note that the degree of $X^n$ is $2^{n+1}=2^{2m}$ and for a generic choice of $Q_i'$,
equation (\ref{eq5}) has no solution of multiplicity bigger than 1.

\qe

By Lemma \ref{lemma2.6}, the map $\rho: X^{2m-1}\to X^{2m-1}$ induces the push forward map
$\rho_*=(-1)^{n}=(-1)^{2m-1}=-1: \Ch_0(X^{2m-1})_{hom}\to \Ch_0(X^{2m-1})_{hom}$. By Lemma \ref{lemma2.3}, $\pi^*\pi_*=\rho_*+1=0:\Ch_0(Y^{2m-1})_{hom}\to \Ch_0(Y^{2m-1})_{hom}$ and  $\Ch_0(Y^{2m-1})_{hom}=0$.

\begin{lemma}
Let $\phi:\widetilde{Y^{2m-1}}\to Y^{2m-1}$ be a resolution of singularity. Then
$\Ch_0(\widetilde{Y^{2m-1})}\cong \Z$.
\end{lemma}
\bp For each singular point $q_i\in Sing(Y^{2m-1})$, $i=1,2,\cdots, 2^{2m}$,  the
exceptional divisor $E_i=\phi^{-1}(q_i)$ has normal crossings
in $\widetilde{Y^{2m-1}}$ and every irreducible component  of $E_i$ is
nonsingular and rational (cf. Corollary after Theorem 1 in \cite{Fujiki}). Since
each singular point in our case is a quotient singularity of
type $\C^{2m-1}/\Z_2$, the exceptional divisor $E_i$ contains exactly
one irreducible component, which is isomorphic to $\P^{2m-2}$ (cf. Remark \ref{remark2.9}).

Set $E=\bigcup_{i=0}^{2^{2m}} E_i$. Since $\Ch_0(E_i)_{hom}=0$ and  that
$E_i$ are mutually disjoint to each other, we get $\Ch_0(E)_{hom}=0$.
Set $U=\widetilde{Y^{2m-1}}-E\cong Y^{2m-1}-\bigcup_{i=0}^{2^{2m}} q_i$.
Then the isomorphism $\Ch_0(\widetilde{Y^{2m-1})}\cong \Z$ follows from the fact that
 $\Ch_0(\widetilde{Y^{2m-1})}_{hom}=0$. This fact  can be seen from
the commutative diagram of Chow groups
$$
\xymatrix{\Ch_0(E)_{hom}\ar[r]\ar[d]^{\cong}& \Ch_0(\widetilde{Y^{2m-1}})_{hom}\ar[r]\ar[d]^{\phi_*}&\Ch_0(U)_{hom}\ar[r] \ar[d]^{\cong}&0\\
\Ch_0(\bigcup_{i=0}^{2^{2m}} q_i)_{hom}\ar[r]&\Ch_0({Y^{2m-1}})_{hom}\ar[r] & \Ch_0(U)_{hom}\ar[r]&0.
}
$$
\qe

\begin{remark}\label{remark2.9}
Each singular point of $Y^{2m-1}$ is the quotient singularity of the same type
as that of  $\C^{2m-1}/\Z_2$, where $\Z_2$ acts on $\C^{2m-1}$
as $$(x_1,x_2,\cdots,x_{2m-1}) \to (-x_1,-x_2,\cdots,-x_{2m-1}).$$
The singular point of $\C^n/\Z_2$ can be resolved by one blow up with
the exceptional divisor $E\cong \P^{n-1}$.

To see this, we first note that all the $\Z_2$-invariant monomials
of $x_1$,  $x_2$, $\cdots$, $x_{n}$ are $x_ix_j$, $1\leq i\leq j\leq n$.
This gives an embedding $X:=\C^n/\Z_2\hookrightarrow \C^N$,
where $N=(^{n+1}_{~2})=\frac{1}{2}n(n+1)$. Let $u_{ij}$, $1\leq i\leq j\leq n$ be
the coordinates of $\C^N$. Then $\C^n/\Z_2$ is the locus of the ideal generated
by all $2\times 2$ minors of the symmetric matrix
$(u_{ij})_{1\leq i,j\leq n}$, where $u_{ij}:=u_{ji}$ if $i>j$.

Let $\widetilde{\C^N}$ be the blow up of $\C^N$ at the origin and
let $\widetilde{X}$ be the proper transform of $X=\C^n/\Z_2$.
A direct calculation shows  that $\widetilde{X}$ is smooth.
The explicit equations for $n=3$ will be given below while the general case is similar.
The exceptional divisor $E$ of $\widetilde{X}\to X$
is just the quadric equation given by those
$2\times 2$ minors in $\P^{N-1}$, i.e, the intersection of $\P^{N-1}$
and $\widetilde{X}$. Note that $E\subset \P^{N-1}$ with the above
defining equations is exactly the image of the  Pl\"{u}cker embedding
$\P^{n-1}\hookrightarrow \P^{N-1}$. Therefore, $E\cong \P^{n-1}$.

Now we write down the details for the case that $n=3$.  In this case
$N=(^4_2)=6$. Let $\widetilde{\C^6}\subset \C^6(u_1,\cdots,u_6)\times \P^5[v_1:\cdots:v_6]$ be defined
by $u_iv_j=u_jv_i$, $1\leq i\neq j\leq 6$. Note that
$X=\C^3/\Z_2\subset \C^6$ is defined by
$$
\left\{\begin{array}{lll}
u_1u_2=u_6^2,\\
u_1u_3=u_5^2,\\
u_2u_3=u_4^2.
\end{array}\right.
$$

Note that $\widetilde{\C^6}$ is covered by affine open sets $(v_i\neq 0)$, $i=1,2,\cdots,
6$.
On the affine open  piece $v_1\neq 0$, $\widetilde{X}$ is defined by
\begin{equation*}
\left\{\begin{array}{lll}
u_i=u_1v_i, i=2,3,\cdots, 6,\\
v_2=v_6^2,\\
v_3=v_5^2,\\
v_2v_3=v_4^2.
\end{array}\right.
\end{equation*}

It is easy to check by the definition of  smoothness that this piece of $\widetilde{X}$ is smooth. Similarly for all other pieces of $\widetilde{X}$.
Therefore $\widetilde{X}$ is smooth.

The exceptional divisor $E$ is defined by the following equations:
\begin{equation*}
\left\{\begin{array}{lll}
u_i=0, i=1,2,\cdots, 6,\\
v_1v_2=v_6^2,\\
v_2v_3=v_5^2,\\
v_2v_3=v_4^2.
\end{array}\right.
\end{equation*}

Hence $E$ is isomorphic to the image of Pl\"{u}cker embedding $\P^2\hookrightarrow \P^5$ and therefore $E\cong \P^2$.

\end{remark}

\section{Application to 1-cycles and codimension two cycles} \label{sec4}
In this section, we deduce a sequence of results on algebraic cycles
and cohomology theories for $Y^n$ as the application of the decomposition
of the diagonal given by Bloch \cite{Bloch},
Bloch and Srinivas \cite{Bloch-S} and the generalization by many others.

First we consider the case  that  $n=2m$ is an  even positive integer.
\begin{corollary}\label{cor1.1}
$\Ch^p(Y^{2m})$ is weakly representable for  $p\leq 2$.
\end{corollary}

\bp It follows from Theorem \ref{Th1.1} and Theorem 1 in \cite{Bloch-S}.

\qe

 The {\bf Hodge Conjecture} for codimension $p$ cycles(denote by $\HC^{q,q}(X)$):
{\sl The rational cycle class map
$$cl^q\otimes{\mathbb{Q}}:{\Ch}^q(X)\otimes{\mathbb{Q}}\rightarrow H^{q,q}(X)\cap
H^{2q}(X,\mathbb{Q})$$ is surjective. }

More generally, let $N^pH^k(X,\Q) \subset H^k(X,\Q)$ be the arithmetic
filtration defined by Grothendieck \cite{Grothendieck} and let
$F^pH^k(X,\C) \subset H^k(X,\C)$ be the Hodge filtration. Set $F^pH^k(X,\Q) : = F^pH^k(X,\C)\cap H^k(X,\Q)$ and denote by
$\widetilde{F}^pH^k(X,\Q)$ the maximal sub-Hodge structure in $F^pH^k(X,\Q)$.
It was shown in \cite{Grothendieck} that $N^pH^k(X,\Q)\subset\widetilde{F}^pH^k(X,\Q)$.

\noindent The {\bf generalized Hodge Conjecture} can be stated as follows (denote by $\GHC(p,k,X)$):
$$
N^pH^k(X,\Q)=\widetilde{F}^pH^k(X,\Q).
$$

\begin{corollary}\label{cor1.2}
The Hodge Conjecture  for $Y^4$ holds.  The generalized Hodge Conjecture
for $\GHC(1,4,Y^{4})$ holds. More generally, the generalized Hodge Conjecture $\GHC(1,2m,Y^{2m})$
for $Y^{2m}$ holds.
\end{corollary}
\bp The first statement follows from Theorem \ref{Th1.1} and Theorem 1 in \cite{Bloch-S}. Similar method can be used to prove $\GHC(1,4,Y^{4})$ and more
general statement $\GHC(1,2m,Y^{2m})$. By Theorem \ref{Th1.1}, $
\Ch_0(Y^{2m})\cong \Z$, we have  $\GHC(1,2m,Y^{2m})$ by
 Corollary 15.23 in \cite{Lewis} or Proposition 5.5 in \cite{Voineagu}.
\qe

\begin{remark}
$\GHC(1,4,Y^{4})$ is the only non-trivial part of the generalized Hodge Conjecture for $Y^4$.
The Hodge Conjecture for 2-cycles and codimension 2 cycles on $Y^{2m}$ holds,
i.e., both $\HC^{2m-2,2m-2}(Y^{2m})$  and $\HC^{2,2}(Y^{2m})$ hold for all
positive integer $m$. However,  both $\HC^{2m-2,2m-2}(Y^{2m})$ and  $\HC^{2,2}(Y^{2m})$
 are trivial if $m>2$ since
both  $H^{4}(Y^{2m})$ and $ H^{4m-4}(Y^{2m})$ are isomorphic to $\Z$.
\end{remark}

Recall that the \textbf{Lawson homology} $L_pH_k(X)$ of $p$-cycles is defined by
$$L_pH_k(X) := \pi_{k-2p}({\cZ}_p(X)) \quad {\rm for}\quad k\geq 2p\geq 0,$$
where ${\mathcal Z}_p(X)$ is provided with a natural topology (cf.
\cite{Friedlander1}, \cite{Lawson1} and \cite{Lawson2}). For general background on
Lawson homology, the reader is referred to \cite{Lawson2}. There
are  natural maps, called \textbf{cycle class maps}
$ \Phi_{p,k}:L_pH_{k}(X)\rightarrow H_{k}(X). $

Define
$$
\begin{array}{lcl}
L_pH_{k}(X)_{hom}&:=&{\rm ker}\{\Phi_{p,k}:L_pH_{k}(X)\rightarrow H_{k}(X)\};\\
L_pH_{k}(X,\Q)_{hom}&:=& L_pH_{k}(X)_{hom}\otimes{\Q};\\
T_pH_{k}(X)&:=&{\rm Image}\{\Phi_{p,k}:L_pH_{k}(X)\rightarrow H_{k}(X)\}; \\
 T_pH_{k}(X,\Q)&:=&T_pH_{k}(X)\otimes\Q.
\end{array}
$$

The \textbf{Griffiths group} of  $p$-cycles is defined to
$$
{\Griff}_p(X):={\Ch}_p(X)_{hom}/{\Ch}_p(X)_{alg},
$$
where ${\Ch}_p(X)_{alg}$ denotes the space of cycles in ${\Ch}_p(X)$
which are algebraically equivalent to zero. Set ${\Griff}^p(X):=\Griff_{\dim X-p}(X)$.
It was shown in \cite{Friedlander1} that $L_pH_{2p}(X)_{hom}\cong \Griff_p(X)$ for any projective variety.

\begin{corollary}\label{cor1.3}
For every positive integer $m$, we have
$L_pH_{k}(Y^{2m},\Q)_{hom}=0$
for $p\leq 1$,  $p\geq 2m-2$ and $k\geq 2p$.
In particular, $\Griff_1(Y^{2m})\otimes \Q=0$ and $\Griff^2(Y^{2m})=0$.
\end{corollary}
\bp  Recall that a theorem of Peters \cite{Peters} says that if
$\Ch_0(Y)\otimes\Q \cong \Q$ for a smooth projective variety $Y$, then
 $L_pH_*(Y)_{hom}=0$ for $p\leq 1$ and all $*$.
Hence the $p\leq 1$ part  follows from  Theorem \ref{Th1.1} and Peters' result.
It was  observed, independently by M. Voineagu \cite{Voineagu} and
the author \cite{Hu}, that Peters' method  could be used to
show $L_pH_*(Y)_{hom}=0$ for $p\geq \dim(Y)-2$ and all $*$ under the
same assumption.
So the  $p\geq 2m-2$ part follows from  Theorem \ref{Th1.1} and the observation.
In particular, $\Griff_1(Y^{2m})\otimes \Q=0$ and $\Griff^2(Y^{2m})\otimes \Q=0$.
Since $\Griff^2(Y^n)$ has no torsion \cite{Bloch-S}, therefore $\Griff^2(Y^{2m})=0$, i.e., homological equivalence and algebraic equivalence coincide for codimension-2 cycles on $Y^{2m}$. The completes the proof of Corollary \ref{cor1.3}.

\qe

As applications of Theorem \ref{Th1.7}, we have similar results for $\widetilde{Y^{2m-1}}$  as those in Corollary \ref{cor1.1},\ref{cor1.2} and \ref{cor1.3}.

\begin{corollary}\label{cor1.8}
$\Ch^p(\widetilde{Y^{2m-1}})$ is weakly representable for  $p\leq 2$.
\end{corollary}

\begin{corollary}\label{cor1.9}
The Generalized Hodge Conjecture
for $\GHC(1,3,\widetilde{Y^{3}})$ holds. More generally, the Generalized Hodge Conjecture $\GHC(1,2m-1,\widetilde{Y^{2m-1}})$
for $Y^{2m}$ holds.
\end{corollary}

\begin{corollary}\label{cor1.10}
For every  integer $m\geq 2$, we have
$L_pH_{k}(\widetilde{Y^{2m-1}},\Q)_{hom}=0$
for $p\leq 1$,  $p\geq 2m-3$ and $k\geq 2p$.
In particular, $\Griff_1(\widetilde{Y^{2m-1}})\otimes \Q=0$ and $\Griff^2(\widetilde{Y^{2m-1}})=0$.
\end{corollary}

Recall that  for $V\subset U$ a Zariski open subset of a quasi-projective
variety $U$, we have the long exact sequence for Lawson homology, i.e.,
\begin{equation}\label{eq6}
\cdots\rightarrow L_pH_k(Z)\rightarrow L_pH_k(U)\rightarrow
L_pH_k(V)\rightarrow L_pH_{k-1}(Z)\rightarrow\cdots
\end{equation}
where $Z=U-V$ (cf. \cite{Lima-Filho}).

By Corollary \ref{cor1.10} and  Equation (\ref{eq6}), we get

\begin{corollary}\label{cor1.11}
For every  integer $m\geq 2$, we have
$L_pH_{k}({Y^{2m-1}},\Q)_{hom}=0$
for $p\leq 1$,  $p\geq 2m-3$ and $k\geq 2p$.
\end{corollary}

\bp Since $Y^{2m-1}$ is a singular variety, the Bloch-Srinivas  method on
decompositions of the diagonal does not work for $Y^{2m-1}$. So we try to
compute $L_pH_{k}({Y^{2m-1}},\Q)_{hom}$ by the localization sequences for Lawson homology.

Set $V:=\widetilde{Y^{2m-1}}-E$ and $E=\cup_{i=1}^{2^{2m}} E_i$. Then
 $V\cong Y^{2m-1}-\cup_{i=1}^{2^{2m}} p_i$, where $p_i$, $i=1,\cdots, 2^{2m}$
 are singular points of $Y^{2m-1}$ and $E_i$,  $i=1,\cdots, 2^{2m}$  are
 the corresponding exceptional divisors. By using Equation (\ref{eq6}) to the
 $U\subset Y^{2m-1}$, we get $L_1H_k(U)\cong L_1H_k(Y^{2m-1})$ and so
 $L_1H_k(U)_{hom}\cong L_1H_k(Y^{2m-1})_{hom}$.

From  the following commutative diagram (cf. \cite{Lima-Filho2}, Prop. 4.9)
{\footnotesize
$$
\xymatrix{L_1H_k(E)\ar[r]\ar[d]^{\cong}& L_1H_k(\widetilde{Y^{2m-1}})\ar[r]\ar[d]^{\Phi_{1,k}}&L_1H_k(U)\ar[r]
\ar[d]^{\Psi_{1,k}}& L_1H_{k-1}(E)\ar[r]\ar[d]^{\cong}&L_1H_k(\widetilde{Y^{2m-1}})\ar[d]^{\Phi_{1,k-1}}\\
H_k(E)\ar[r]& H_k(\widetilde{{Y^{2m-1}}})\ar[r] & H_k^{BM}(U)\ar[r]&H_{k-1}(E)\ar[r]&H_k(\widetilde{Y^{2m-1}})}
$$
}
where $H_k^{BM}(U)$ is the Borel-Moore homology of $U$, and the injectivity
of $\Phi_{1,k}\otimes \Q$ (i.e., $L_1H_{k}(\widetilde{Y^{2m-1}},\Q)_{hom}=0$
by Corollary \ref{cor1.10}), we get the injectivity of
$\Psi_{1,k}\otimes \Q$ by the Five Lemma.
\qe

\section{Low dimensional examples}
For a smooth complex projective variety, we set $h^{i,j}(X):=\dim_{\C} H^{i,j}(X)$.

The case $n=1$ is trivial. In this case $Y^1\cong \P^1$.
In the case $n=2$, all $Y^2$ are  Enrique surfaces.
It was proved in \cite{Bloch-K-L} that all Enrique surfaces $S$ satisfy $\Ch_0(S)\cong \Z$.

The next case is $n=3$. In this case, $X^3$ (for simplicity denote by $X$ in
this paragraph) is the complete intersection of
4 quadric hypersurfaces in $\P^7$. By the adjunction formula, the canonical
bundle $K_X$ of $X$ is trivial and so $H^{3,0}(X)\cong \C$ and
$h^{3,0}(X)=1$. The Euler class $\chi(X)$ of $X$ is the top
Chern class of $X$ (Gauss-Bonnet Theorem). Let $h$ be the hyperplane class
of $X$. The total Chern class
$c(X):=1+c_1(X)+c_2(X)+c_3(X)=(1+h)^8(1+2h)^{-4}|_X$ and so
$c_3(X)=-8h^3=-128$ since $h^3|_X=\deg(X)=16$.  Hence $\chi(X)=-128$.
This together with Lefschetz hyperplane theorem implies $b_3(X)=132$.
By the Hodge decomposition of $H^3(X,\C)$,
we get $h^{3,0}(X) + h^{2,1}(X) + h^{1,2}(X) + h^{0,3}(X) = 132$.
So $h^{2,1}(X)=h^{1,2}(X)=65$.

Since $\pi:X^3\to Y^3$ is an involution with 16
isolated fixed point, we have
$\chi(X^3)-16=2(\chi(Y^3)-16)$ and so $\chi(Y^3)=-56$. Note that
$\widetilde{Y^3}\to Y^3$ is the resolution of singularity with exceptional
divisor $E=\cup_{i=1}^{16} E_i$, where $E_i\cong \P^2$ since it is the
exceptional divisor  of the resolution of singularity on $\C^3/\Z_2$
(cf. Remark \ref{remark2.9}). So
$\chi(\widetilde{Y^3})-16\chi(\P^2)= \chi(Y^3)-16$, i.e.,
$\chi(\widetilde{Y^3}) = -24$. Since $b_2(Y^3)=1$, we get
$b_2(\widetilde{Y^3}) = 17$. Hence $b_1(\widetilde{Y^3})=0$,
we get $b_3(\widetilde{Y^3})=60$. We get $h^{1,2}(\widetilde{Y^3}) = h^{2,1}(\widetilde{Y^3})=30$ since $h^{3,0}(\widetilde{Y^3})=0$.

Recall that Suslin's Conjecture on Lawson homology states that:
\emph{For any abelian group $A$ and smooth quasi-projective variety $X$ of dimension $n$, the map $L_pH_k(X,A) \to H^{BM}_k(X,A)$ is an isomorphism for $k\geq n+p$ and a monomorphism for $k=n+p-1$.
} Here $H^{BM}_k(X,A)$ means Borel-Moore homology with coefficient in $A$.

By the above computation, Theorem \ref{Th1.7} and  the
Proposition 5.1 in \cite{Voineagu},  we have
\begin{corollary}
 For $\widetilde{Y^3}$ above,  we have the following statements:
 \begin{enumerate}
\item $L_1H_2(\widetilde{Y^3})\cong H_2(\widetilde{Y^3}) \cong \Z^{17}$.
 \item  $L_1H_3(\widetilde{Y^3})\cong H_3(\widetilde{Y^3})\cong \Z^{60}$.
 \item $L_2H_4(\widetilde{Y^3})\cong L_1H_4(\widetilde{Y^3})\cong H_4(\widetilde{Y^3}) \cong \Z^{17}$.
  \item    $L_3H_6(\widetilde{Y^3})\cong L_2H_6(\widetilde{Y^3})\cong L_1H_6(\widetilde{Y^3})\cong H_6(\widetilde{Y^3}) \cong \Z$.

  \item All other $L_pH_*(\widetilde{Y^3})$ ($p\geq 1$) are trivial.

  In particular,  Suslin's Conjecture for $\widetilde{Y^3}$ holds.
 \end{enumerate}
\end{corollary}

Next case is $n=4$. In this case, it can be calculated that
$h^{4,0}(X^4)=1$, $h^{3,1}(X^4)=h^{1,3}(X^4)=151$ and $h^{2,2}(X^4)=652$
since $X^4$ is a complete intersection. Since $h^{4,0}(Y^4)=0$, we get
 $h^{3,1}(Y^4)=h^{1,3}(Y^4)=75$ and $h^{2,2}(Y^4)=326$ by Riemann-Roch-Hirzebruch
 theorem for orbit spaces (cf. \cite{Atiyah-Singer}, 4.7),
Since $h^{3,1}(Y^4)\neq 0$, $\Ch_1(Y^4)$ is not weakly representable.
In particular, $\Ch_1(Y^4)_{hom}\otimes\Q$ is nontrivial. From this and the proof
of Theorem \ref{Th1.1}, we obtain the Chow group of 1-cycles
for  $X^4/\langle\tau_i\rangle$ is not weakly representable,
although $X^4/\langle\tau_i\rangle$ is a rationally connected variety for
each $i=0,1,\cdots, n$.

In this case we can say a little more about Lawson homology of $Y^4$.   By
By the above computation,  Theorem \ref{Th1.1} and  the
Proposition 5.3 in \cite{Voineagu},  we have the following result.
\begin{corollary}For $Y^4$ above, we have the following statements
\begin{enumerate}
\item $L_1H_2(Y^4)_{\mathbb{Q}}{\cong}H_2(Y^4)_{\mathbb{Q}}\cong\Q.$


\item $L_2H_4(Y^4)\hookrightarrow L_1H_4(Y^4)\cong H_4(Y^4)\cong\Z^{476}.$


\item $L_3H_6(Y^4)\cong L_2H_6(Y^4)\cong L_1H_6(Y^4)\cong H_6(Y^4)\cong\Z.$


\item $L_4H_8(Y^4)\cong L_3H_8(Y^4)\cong L_2H_8(Y^4)\cong L_1H_8(Y^4)\cong H_8(Y^4)\cong\Z.$

\item All other $L_pH_*(Y^4)$ ($p\geq 1$) are trivial.
\end{enumerate}

In particular,  Suslin's Conjecture for $Y^4$ holds.
\end{corollary}

\begin{remark}
One can show that $H^{1,n-1}(Y^{2m},\Q)$ is nontrivial for all  $n=2m\geq 4$. Therefore $\Ch_1(Y^{2m})_{hom}\otimes\Q\neq 0$ and Theorem \ref{Th1.1} is the best result we could obtain. Similarly,
One can show that $H^{1,n-1}(Y^{2m-1},\Q)$ is nontrivial for all  $n=2m-1\geq 3$
and Theorem \ref{Th1.7} is the best result.

\end{remark}

Department of Mathematics, Massachusetts Institute of Technology, Room 2-363B,
 77 Massachusetts Avenue,
 Cambridge, MA 02139, USA

Email: {wenchuan@math.mit.edu}

\begin{thebibliography}{AAAA}
\bibitem[AS]{Atiyah-Singer}
 M. F. Atiyah and I. M. Singer,
{\sl The index of elliptic operators. III.}
Ann. of Math. (2) 87 1968 546--604.

\bibitem[B1]{Bloch}
S. Bloch,
{\sl Lectures on algebraic cycles.}
Duke University Mathematics Series, IV. Duke University, Mathematics Department, Durham, N.C., 1980. 182 pp. (not consecutively paged).

\bibitem[B2]{Bloch2}
S. Bloch,
{\sl Some elementary theorems about algebraic cycles on Abelian varieties.}
Invent. Math. 37 (1976), no. 3, 215--228.

\bibitem[BKL]{Bloch-K-L}
S. Bloch, A. Kas, and D. Lieberman,
{\sl Zero cycles on surfaces with $p\sb{g}=0$.}
Compositio Math. 33 (1976), no. 2, 135--145.


\bibitem[BS]{Bloch-S}
S. Bloch and V. Srinivas,
{\sl Remarks on correspondences and algebraic cycles.}
Amer. J. Math. 105 (1983), no. 5, 1235--1253.

\bibitem[C]{Campana}
F. Campana,
{\sl Connexit\'e rationnelle des vari\'et\'es de {F}ano.}
Ann. Sci. \'Ecole Norm. Sup. (4) 25 (1992), no. 5, 539--545.

\bibitem[Fr]{Friedlander1} E. Friedlander, {\sl Algebraic cycles, Chow
varieties, and Lawson homology.}  Compositio Math. 77 (1991), no. 1,
55--93.

\bibitem[Fuj]{Fujiki}
A. Fujiki,
{\sl On resolutions of cyclic quotient singularities.}
Publ. Res. Inst. Math. Sci. 10 (1974/75), no. 1, 293--328.

\bibitem[Ful]{Fulton}
W. Fulton, {\sl Intersection theory.} Second edition, Springer-Verlag,
Berlin, 1998.


\bibitem[G]{Grothendieck} A. Grothendieck,
{\sl Hodge's general conjecture is false for trivial reasons.} Topology 8 1969 299--303.

\bibitem[H]{Hu}
W. Hu, {\sl A note on Lawson homology for smooth varieties with small Chow groups.}
 arxiv:math/0602516



\bibitem[IM]{Inose-Mizukami}
H.Inose and  M. Mizukami,
{\sl Rational equivalence of $0$-cycles on some surfaces of general type with $p\sb{g}=0$}.
Math. Ann. 244 (1979), no. 3, 205--217.


\bibitem[KMM]{Kollar-Miyaoka-Mori}
J. Koll\'{a}r, Y. Miyaoka and S. Mori,
{\sl Rational connectedness and boundedness of Fano manifolds.}
J. Differential Geom. 36 (1992), no. 3, 765--779.

\bibitem[L1]{Lawson1}
B. Lawson, {\sl Algebraic cycles and homotopy theory.}, Ann. of
Math. {\bf 129}(1989), 253-291.

\bibitem[L2]{Lawson2}B. Lawson, {\sl Spaces of algebraic
cycles.} pp. 137-213 in Surveys in Differential Geometry, 1995
vol.2, International Press, 1995.

\bibitem[Le]{Lewis}
J. D.  Lewis,
{\sl A survey of the Hodge conjecture.} (English summary)
Second edition. Appendix B by B. Brent Gordon. CRM Monograph Series, 10.
American Mathematical Society, Providence, RI, 1999. xvi+368 pp. ISBN: 0-8218-0568-1

\bibitem[LF]{Lima-Filho}P. Lima-Filho, {\sl Lawson homology for quasiprojective
varieties.} Compositio Math.  84  (1992),  no. 1, 1--23.

\bibitem[LF2]{Lima-Filho2}
P. Lima-Filho,
{\sl On the generalized cycle map.} (English summary)
J. Differential Geom. 38 (1993), no. 1, 105--129.


\bibitem[M]{Mumford}
D. Mumford,
{\sl Rational equivalence of $0$-cycles on surfaces.}
J. Math. Kyoto Univ. 9 1968 195--204.
\bibitem[Pe]{Peters} C. Peters, {\sl  Lawson homology for varieties with
small Chow groups and the induced filtration on the Griffiths
groups.}  Math. Z. 234 (2000), no. 2, 209--223.

\bibitem[PV]{Paranjape-Srinivas}
K. H. Paranjape and V. Srinivas,
Algebraic cycles. Current trends in mathematics and physics, 71--86, Narosa, New Delhi, 1995.

\bibitem[R]{Roitman}
A.A. Ro\u\i tman,
{\sl Rational equivalence of zero-dimensional cycles.} (Russian)
Mat. Sb. (N.S.) 89(131) (1972), 569--585, 671.

\bibitem[Vo]{Voineagu} Mircea Voineagu,
{\sl Semi-topological K-theory for certain projective varieties.}
Preprint. arxiv.org/abs/math/0601008


\bibitem[Vs]{Voisin}
C. Voisin,
{\sl Sur les z\'ero-cycles de certaines hypersurfaces munies d'un
automorphisme }
Ann. Scuola Norm. Sup. Pisa Cl. Sci. (4) 19 (1992), no. 4, 473--492.

\end{thebibliography}
\end{document}